\documentclass[11pt,letterpaper]{amsart}

\usepackage{amssymb} 
\usepackage{amsthm} 
\usepackage{amsfonts}
\usepackage{amsmath}
\usepackage{amstext}
\usepackage[dvips]{graphics}
\usepackage{epsfig}

\numberwithin{equation}{section}

\theoremstyle{plain}%default
\newtheorem{thm}{Theorem}[section] 

\newtheorem{cor}[thm]{Corollary}
\newtheorem{lem}[thm]{Lemma}

\newtheorem{theorem*}{Theorem}[]

\theoremstyle{definition}
\newtheorem{defn}[thm]{Definition}
\newtheorem{conv}[thm]{Convention}

\newtheorem{example}[thm]{Example}

\theoremstyle{remark}
\newtheorem{rem}[thm]{Remark}

\newcommand{\C}{\mathbb{C}}

\newcommand{\inv}{^{-1}}

\newcommand{\val}{\nu}
\let\mathscr\mathcal
\makeatletter
\@input{mathrsfs.sty}
\makeatother

\newcommand{\mer}{\mathscr M}
\newcommand{\pmer}{\mathcal P}

\newcommand{\post}{\rightarrow}
\newcommand{\allpost}{\mathcal A}
\newcommand{\rep}{\mathcal R} 
\newcommand{\T}{\mathcal T}

\title{Newton-Puiseux Roots of Jacobian Determinants}

\author{Tzee-Char Kuo}
\thanks{The first author is partially supported by an ARC Large Grant.}
\address{School of Mathematics and Statistics, University of Sydney,
  Sydney, NSW, 2006, Australia }
\email{tck@maths.usyd.edu.au; fax:61-2-93514534.}
\author{Adam Parusi\'nski}
\address {D\'epartement de Math\'ematiques, U.M.R. 6093 du C.N.R.S,
Universit\'e d'Angers, 2, bd Lavoisier, 49045 Angers Cedex, France}
\email{parus@tonton.univ-angers.fr}


\begin{document}

\begin{abstract} 
Let $f(x,y), g(x,y)$ denote either a pair of holomorphic 
function germs, or a pair of monic polynomials in $x$ whose
coefficients are Laurent series in $y$. A {\it polar root} is
a Newton-Puiseux root, $x=\gamma(y)$, of the Jacobian
$J=f_yg_x-f_xg_y$, but not a root of $f\cdot g$.

We define the tree-model, $T(f,g)$, for the pair, 
using the set of contact orders of the Newton-Puiseux roots of $f$ and
$g$. Our main results (\S \ref{Results}) describe how the $\gamma$'s
climb, and leave, the tree (like vines).
We also show by two examples (\S \ref{C}) that when the tree has what
we call collinear points or bars, the way the $\gamma$'s leave the
tree is not an invariant of the tree; this phenomenon is in sharp
contrast to that in
the one function case where the tree $T(f)$ completely determines how
the polar roots split away (\cite {kuo-lu}, \cite {kuo-parusinski}).

Our results yield a factorisation of the Jacobian determinant in $\C\{x,y\}$ 
(\S \ref{F}). As in the one function case, the factors need not be invariants,
 nor irreducible. However, some factors do yield invariant truncations and 
intersection multiplicities (\S \ref{Inv}). 
\end{abstract}

\maketitle

\vspace{ 0.8 truecm}

Take two holomorphic germs $f, g: (\C^2,O) \post (\C,O)$, and a coordinate 
system $(x,y)$.  The Newton-Puiseux factorisations are of the form
\begin{equation}\label{fandg}\begin{split}  
f(x,y)=u(x,y)\cdot y^{E_1}\cdot \prod_{i=1}^p [x-\alpha_i (y)],\quad
E_1\ge 0,\\
g(x,y)=u'(x,y)\cdot y^{E_2}\cdot \prod_{j=1}^q [x-\beta_j (y)],\quad
  E_2\ge 0,
\end{split}\end{equation}
where $u$, $u'$ are units, $\alpha_i$, $\beta_j$ are fractional power series
with $O_y(\alpha_i)>0$, $O_y(\beta_j)>0$.

We shall also write $\alpha_1, \ldots, \alpha_p, \beta_1,\ldots,
\beta_q$ as $\lambda_1,\ldots, \lambda_N$, $N:=p+q$.

\begin{defn}%\label
{\it A polar root} of the pair $(f,g)$, relative to the coordinate
system $(x,y)$, is a
Newton-Puiseux root, $x=\gamma(y)$, $O_y(\gamma)>0$, 
of the Jacobian determinant 
\begin{equation*}
J(x,y):=J_{(f,g)}(x,y) := \biggl| \begin{matrix}
f_y  & f_x  \\
g_y & g_x  \end{matrix} \biggr| ,  
\end{equation*}
which is not one of the $\lambda_k$'s,  that is:
\begin{equation*}
J(\gamma(y),y)=0, \quad f(\gamma(y),y) g(\gamma(y),y)\ne 0.
\end{equation*}
\end{defn}
\medskip
Polar curves play an important r\^ole in Singularity Theory;  
they have been intensively studied by many authors from different 
perspectives. See, e.g.\,\cite{assi}, \cite{teissier}, \cite{Bo}, 
\cite{kuo-parusinski}, \cite {le}. 

It is easy to see that the Newton-Puiseux roots of $J(x,y)$ are
the polar roots plus the multiple roots of the product function 
$f(x,y)g(x,y)$. 

In this paper we shall assume that $f(x,y)g(x,y)$ 
has only simple roots: $\lambda_i\ne \lambda_j$ if $i\ne j$.

We shall use the contact orders 
$O(\lambda_s,\lambda_t):=O_y(\lambda_s(y)-\lambda_t(y))$, $1\le s, t \le N$, 
to associate to the pair $(f,g)$ a combinatorial object: the tree model 
$T(f,g)$.  We use the tree model to analyse the contact orders between the 
polar roots and the roots $\lambda_i$, so that we can visualize how the polar 
roots climb (like vines) along the tree and how they leave it.  
Each bar of the tree gives rise to a rational function $\mer(z)$ of one 
complex variable. This function depends not only on the orders of
contact between the $\lambda_k's$ but also on the coefficients.  
In particular, its zeros play an important r{\^o}le since we can identify 
them with some of the places where the polar roots leave the tree.  But it may 
happen that this function $\mer (z)$ is identically zero. 
 Then the bar is called {\it collinear} and, as we show in \S \ref{C}, it 
is impossible, in general, to know precisely how the polar roots
climb this bar and leave the tree.  
 
\smallskip
We have divided our results into two parts. Part 1, in \S
\ref{Results}, contains three theorems (Theorems T, N, and C) and their
corollaries. They describe the positions of the polar roots relative to
$T(f,g)$. The proofs depend heavily on the classical Theorem of Rouch\'e. 
We show, in particular, that if a polar root leaves the tree on a non-collinear
 bar, it must do so at a ``pure'' zero of the associated rational 
function (see the next section for this terminology).  This will be
used in Part 2 to obtain a factorisation of the Jacobian determinant.

Part 2, in \S \ref{F} and \S \ref{Inv}, contains results which describe how
$J(x,y)$ can be factored in $\C \{x,y\}$, and how to compute the
intersection multiplicities of each factor with the zero sets of $f$, $g$, 
and $f\cdot g$. The factors can be reducible. Amongst the objects we come 
across, we shall carefully distinguish those which are invariants of the tree 
from those which are not (\S \ref{Inv}). 

\smallskip
\emph{Our results generalise that in the one function case. }
More specifically, by taking
$g(x,y)=y$, $J(x,y)$ reduces to $f_x$, and $T(f,y)=T(f)$, the
tree-model defined in \cite{kuo-lu}. 
The curve $f_x=0$ is called a polar curve, it has been studied since 
the time of M. Noether. The components, defined
by the irreducible factors of $f_x$ in $\C\{x,y\}$, are called the
polar branches. 
In \cite{Pham}, Pham showed that the Zariski equisingularity type 
(\cite{Zariski}) of the polar curve, and that of the polar branches, need not 
be determined by that of $f=0$. 
(If $f$ is generic in its equisingularity class then the 
equisingularity class of the generic polar curve is described in 
\cite{casas}.)

However,
for the contact orders of the polar roots with
the roots of $f=0$, the story is different.  The set of contact orders,
$C(f,f_x):=\{O(\alpha_i,\gamma_j)\}$, between the roots of $f$ and
that of $f_x$ can be calculated using the tree-model $T(f)$ {\it
  alone}. This is proved in \cite{kuo-lu},
and, for irreducible $f$, also in \cite{merle}.

\emph{ Therefore}, the contact order set $C(f,f_x)$ \emph{is an  invariant  
of the equisingularity type of $f$}.

(Attention should be paid to the rather subtle distinction between
polar roots, which are studied in this paper, and the well-established
notion of polar branches; the former are fractional power series, the
latter are primes in $\C \{x,y\}$ generated by the former.
This difference between polar roots
and polar branches has eluded some experts.)

A more detailed account of the one function case 
will be given in Section \ref{Onefunction}. 

The main results of this paper have been announced in \cite{kuo-parusinski2}.

\smallskip
\emph{In the general case,} $T(f,g)$ may have collinear points and bars 
(in the one function case, all bars of $T(f)$ are purely non-collinear), 
and then   we encounter a completely new phenomenon.  Namely, 
it may not be possible anymore to know precisely where some of the polar roots 
leave the tree. That is, the set of contact orders
$\{O(\lambda_i,\gamma_j)\}$ between the roots $\lambda_i$ of $f\cdot g$ and the
polar roots $\gamma_j$ need not be determined by $T(f,g)$. 
We give examples in \S \ref{C}.

There are many other tree-models in the one function case (Cassas-Alvero,
Eggars, Wall, etc). The models $T(f)$,
$T(f,g)$, we use here can better express how the polar roots split away from
the tree. Our definitions are simple
but do not use the language of Graph Theory.

\smallskip
\noindent
{\bf Conventions: }
A fractional power series $\lambda(y)$ will also be called an "arc".  
If $O(\lambda,\mu)>q$, we write $\lambda\equiv \mu \mod q^+$.
We use $O(y^+)$ to represent a quantity which, as $y\to 0$, has
the same order as $y^e$, for some $e>0$.
Finally, "$+\cdots$" will mean "plus higher order terms".

%%%%%%%%%%%%%%%%%%%%%%%%%%%%%%%%%%%%%%%%%%%%%%%%%%%%%%%%%%%%%%%%%
\medskip
\section{The tree model $T(f,g)$.}\label{Treemodel}
\medskip 

The tree model is a geometric object that allows us to visualize the 
numerical data given by the contact orders 
$O(\lambda_s,\lambda_t)$ between the roots of $f\cdot g$ and then 
between a given arc $x=\xi (y)$ and the $\lambda_i$'s. 
The construction of $T(f,g)$ is as follows (compare
\cite{kuo-lu}). First, draw a horizontal bar, denoted by $B_*$, and
call it the {\it ground bar} (the soil). Then draw a vertical line
segment on $B_*$ as the {\it main trunk} of the tree.  Mark $[p,q]$ alongside
the trunk to indicate that $p$ $\alpha_i$'s and $q$
$\beta_j$'s are bundled together.

Let $h_0:= \min \{O(\lambda_i,\lambda_j) | 1\le i,j\le N\}$. Then draw a bar, 
$B_0$, on top of the main trunk. Call $h(B_0):= h_0$ the
{\it height} of $B_0$. We define $h(B_*):=0$.

The roots $\lambda_k, 1\le k\le N$, are divided into equivalence
classes modulo $h_0^+$.  We then represent each equivalence class
by a vertical line segment drawn on top of $B_0$.  Each is called
a {\it trunk}.  

If a trunk consists of $s$ $\alpha_i$'s and $t$ $\beta_j$'s
($s\ge 0$, $t\ge 0, s+t\ge 1$), we say it has {\it bimultiplicity} $[s,t]$, 
and mark $[s,t]$ alongside. We call $s+t$ the {\it total multiplicity}.

Now, the same construction is repeated recursively on each trunk,
getting more bars, then more trunks, etc..
The height of each bar, the bimultiplicity and the total multiplicity of each 
trunk, are defined likewise.

The construction terminates at the stage where the bars have infinite
height.  We shall omit drawing bars of infinite height.

\begin{example}\label{ex1.1}
Take constants $A\ne 0\ne B$, integers $0<e<E$. Then consider 
\begin{eqnarray*}
f(x,y) = (x+y)(x-y^{e+1} +Ay^{E+1}) (x+y^{e+1} +By^{E+1}), \\
g(x,y)= (x-y) (x-y^{e+1} -Ay^{E+1}) (x+y^{e+1} -By^{E+1}) .
\end{eqnarray*}
The tree model $T(f,g)$ is shown in Fig.1 with $h(B_0)=1$, 
$h(B_1)= e+1$, $h(B_2)=h(B_3) = E+1$. There are six roots $\lambda_k$, 
hence six bars of infinite height.  
The notations "$\circ$" and "$\times$" will be defined in Convention 
\ref{lau}. \end{example}

\epsfxsize=8cm
\epsfysize=5cm
$$\epsfbox{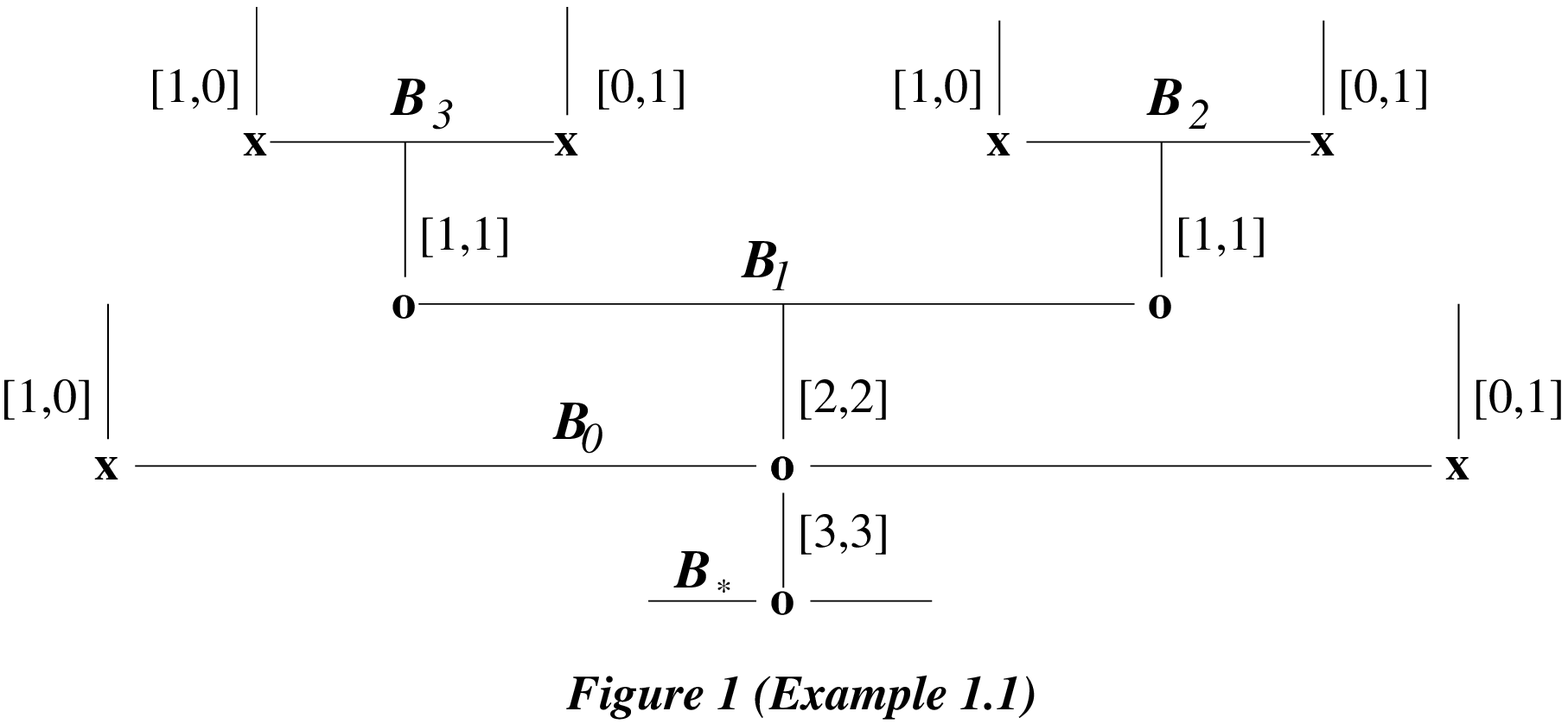}$$

Tracing upward from the main trunk to a bar of infinite height
amounts to identifying a root $\lambda_k$. The heights of the
bars coming  across on the way up are the contact orders of $\lambda_k$ with 
the other roots.
\medskip

Take a bar $B$, with finite height $h:=h(B)$.  Take a root $\lambda_k$
whose modulo $h^+$ class is a trunk on $B$.  Let $\lambda_B(y)$
denote $\lambda _k(y)$ with all terms $y^e$, $e\ge h$, omitted. (In
particular, $\lambda_{B_*}(y)=0$.)
Clearly, $\lambda_B$ depends only on $B$, not on the choice of
$\lambda_k$.  We can then write
\begin{equation*}
\lambda_{k}(y) = \lambda_B(y) + c y^{h(B)} +  \cdots ,
\qquad c\in \C, 
\end{equation*} 
where $c$ is {\it uniquely} determined by $\lambda_k$.  

We say the trunk $T$ which contains $\lambda_k$ {\it grows} on $B$ {\it at} 
$c$.  If $T$ has bimultiplicity $[s,t]$, we also say $c$ has \emph{
bimultiplicity} $[s,t]$ on $B$. The main trunk grows on $B_*$ at $0$.

Let $B^*$ be the bar on top of $T$. 
As in \cite{IKK}, we say $B^*$ is a {\it postbar} of $B$, {\it supported at} 
$c$,
and write: $B\perp _c B^*$. In case $c$ need not be specified, we simply 
write: $B\perp B^*$.

We say $B'$ {\it lies above} $B$ if there is a postbar
sequence:\,$B\perp B_1\perp \cdots \perp B_q \perp B'$; in this case,
if $B\perp_c B_1$, we also say $B'$ {\it lies over} $c$, or $c$ {\it lies
below} $B'$.
\begin{defn}
Take any arc $\xi$.  If $\xi$ has the form
\begin{equation*}
\xi(y) = \lambda_B(y) + a y^{h(B)} +  \cdots ,
\quad a\in \C, 
\end{equation*}
then we say $\xi$ {\it climbs} over $B$ {\it at} $a$ (like a vine). In
this case, if no trunk grows at $a$ we say $\xi$ \emph{leaves} the tree
on $B$ at $a$.
 
If $O(\xi,\lambda_B)<h(B)$, we say $\xi$ is \emph{bounded} by $B$.  
\end{defn}

\begin{conv}\label{copy}
Take a bar $B$. We shall identify $z\in \C$ with the arc 
$\lambda_B(y)+zy^{h(B)}$; and shall also use $B$ to denote the set of
such arcs. (Intuitively, $B$ is a copy of $\C$.)
\end{conv}

%%%%%%%%%%%%%%%%%%%%%%%%%%%%%%%%%%%%%%%%%%%%%%%%%%%%%%%%%%%%%%%
\bigskip
\section{The main results.}\label{Results}
\medskip 

In this section we describe the possible positions of a polar root of 
the pair $(f,g)$ with 
respect to the tree $T(f,g)$.  
Take a bar $B$, $h(B)<\infty$. Take a germ $F(x,y)$ and a generic $z\in \C$. 
Let 
\begin{equation*}
\val_F(B) := O_y(F(\lambda_B(y)+ z y^{h(B)},y)),
\end{equation*}
and, for any fractional power series $\eta(y)$,
\begin{equation*}
\val_F(\eta) := O_y(F(\eta(y),y)).
\end{equation*}
In particular, $\nu_f(B_*)=E_1$, $\nu_g(B_*)=E_2$, by \eqref{fandg}.
\medskip

Let $T_k$, $1\le k\le l$, be the set of  trunks on $B$.  Suppose
$T_k$ grows at $z_k$, having bimultiplicity $[p_k,q_k]$.  
We write
\begin{equation*}
\Delta _B(z_k):= 
\biggl| \begin{matrix}
\val_f(B)             & p_k  \\
\val_g(B)   & q_k  
\end{matrix} \biggr|,\quad 1\le k\le l,  
\end{equation*}
and call 
\begin{equation*}
\mer_B(z) := \sum_{k=1}^l \frac {\Delta _B(z_k)}{z-z_k}, \quad
z\in \C, 
\end{equation*}
the {\it rational function associated to $B$}.  
(Those terms with $\Delta_B(z_k)=0$ can be omitted.)  

\begin{defn}
We say $z_k$, $1\le k\le l$, is a {\it collinear} point on $B$ if
$\Delta_B(z_k)=0$, otherwise, $z_k$ is called 
{\it non-collinear}. 
\end{defn}

Let $C(B)$ and $N(B)$ denote respectively the set of collinear and 
non-collinear points:
\begin{equation}\label{allzk}
C(B)\cup N(B)= \{z_1,\ldots, z_l\} .
\end{equation}
Their (finite) cardinal numbers are denoted by $c(B)$ and $n(B)$ respectively.

\medskip
\begin{conv}\label{lau}
A collinear point will be indicated by
$\circ$; a non-collinear point by $\times$.  
\end{conv}

\begin{defn}\label{Col}
We call $B$ a \emph{ collinear bar} if $\Delta_B(z_k)=0$ for all $k$, 
$1\le k \le l$.  Otherwise we call $B$ \emph{non-collinear}.  We say $B$ is 
\emph{ purely non-collinear} if $C(B)=\emptyset \ne N(B)$. 
\end{defn}

In Example \ref{ex1.1}, $B_1$ and $B_*$ are collinear, $B_2,B_3$ are purely
non-collinear. 

If $\mer_B(z)=0$, we say $z$ is a \emph{mero-zero} on $B$. Let $m_B(z)$ denote
 its multiplicity.
Let $M(B)$ denote the set of mero-zeros on $B$. We write:
\begin{equation*}
m(B):=\sum_{z\in M(B)} m_B(z).
\end{equation*}

Suppose $N(B) \ne \emptyset$. A non-collinear $z_k$ is a simple pole, hence
not a mero-zero:
\begin{equation}
N(B)\cap M(B)= \emptyset, \quad n(B) \ge m(B) +1  .
\end{equation}
On the other hand it may happen that $C(B)\cap M(B) \ne \emptyset$.  
If $z \in M(B) \setminus C(B)$, we say $z$ is a \emph{pure mero-zero}.

It can happen that $\Delta(z_k)=0$ for all $k$.  In this case
$\mer_B\equiv 0$, 
$N(B)=\emptyset$, and $M(B) = \C$.

It can also happen that $M(B) = \emptyset$.
For example, for $f(x,y) =x$ and $g(x,y) =x^2-y^2$, there is only
one bar $B$ of height $1$ and we have:
\begin{equation*}
\mer_B(z) = \frac {
\bigl| \begin{smallmatrix}
1   & 1  \\
2 & 0
\end{smallmatrix} \bigr|} {z} +
\frac {
\bigl| \begin{smallmatrix}
1 & 0  \\
2 & 1 
\end{smallmatrix} \bigr|} {z-1} +
\frac{ \bigl| \begin{smallmatrix}
1 & 0  \\
2 & 1 
\end{smallmatrix} \bigr|} {z+1}= \frac 2 {z(z^2-1)} .
\end{equation*} 

Take a non-collinear bar $B$. 
We define formally the {\it total multiplicity function} by
\[
\tau_B(z)=
\begin{cases}
p_k+q_k,&\text{if $z=z_k;$}\\
0,&\text {otherwise,}
\end{cases} 
\]
and the {\it mero-multiplicity function} by
\[
\mu_B(z)=
\begin{cases}
m_B(z),&\text{if}\;  z\in M(B);\\
-1, &\text{if}\;  z\in N(B);\\
0, &\text{otherwise}.
\end{cases}
\]  

We also write
\begin{equation}
\tau(B):=\sum_{z\in \C} \tau_B(z), \qquad \mu(B)
:=\sum_{z\in \C}\mu_B(z). 
\end{equation}

Note that, obviously,
\begin{equation}
\mu(B)= m(B)-n(B),\quad \text{a negative integer.}
\end{equation}

Let $\T_B(z)$ denote the total number of polar roots 
(counting multiplicities) which climb over $B$ at $z$, and let
$\T(B)$ denote that of those which climb over $B$.

\medskip
\noindent{\bf Theorem T. }
{\em  
Let $B$ be a non-collinear bar. Then
\begin{equation}
\T_B(z)=\tau_B(z)+\mu_B(z), \quad z\in \C,
\end{equation}
and, consequently,
\begin{equation}
\T(B)=\tau(B)+\mu(B).
\end{equation}}

\begin{cor}\label{P}
If a polar root climbs over $B$ at $z$, then $z\in N(B)\cup C(B)\cup M(B).$
\end{cor}

\begin{cor}\label{M}
Suppose $z$ is a pure mero-zero on $B$. Then there are exactly $m_B(z)$ polar 
roots (counting multiplicities) climbing over $B$ at $z$. (Thus, they all leave
 the tree at $z$.)
\end{cor}

\begin{cor}\label{TL}
Suppose $\sum _{z_k\in N(B)} \Delta_B(z_k)\ne 0$.  
Then $m(B) + 1 = n(B)$
and $$\T (B) = \sum_{k=1}^l (p_k + q_k) - 1.$$  
In particular, if $pE_2-qE_1\ne 0$ then the total number of polar roots is 
$p+q-1$.
\end{cor}

\begin{cor}\label{ST} 
Suppose that $B$ is on the top of trunk $T$ of bimultiplicity $[s,0]$ 
(resp.  $[0,t]$) and $\val_g(B) \ne 0$ (resp. $\val_f(B) \ne 0$).
Then 
$B$ is purely non-collinear, $m(B) +1 = n(B)$, and $\T (B) = s-1$ 
(resp. $t-1$).  
\end{cor}

Now we shall study how the polar roots climb the tree, and where they
leave.
% (When there are collinear bars, this is not entirely determined
%by the tree, see \S \ref {C}.)  
The simplest cases have already been dealt with in the above
corollaries. (See also \S \ref{C}.)

\medskip
\noindent {\bf Theorem N. }
{\em Take $z\in N(B)$.  Let $B^*$ be the postbar of $B$ supported at $z$. Then 
\begin{equation}\label{m+1=n}
m(B^*) +1 = n(B^*) .
\end{equation}
In particular, $B^*$ is non-collinear. Moreover, 
every polar root which climbs over $B$ at $z$ must
also climb over $B^*$.  That is, there is no polar
root, $\gamma$, such that
\begin{equation}\label{between}
h(B) < O(\gamma, \lambda_{B^*})< h(B^*) .
\end{equation}  }

\begin{defn} 
Take $c\in C(B)$. A set of non-collinear bars 
$\{ \bar B_1, \cdots, \bar B_l \}$ is called a 
(non-collinear) {\it cover} of $c$ if the following holds:
\begin{itemize} 
\item [(i)]
Each $\bar B_s$ lies over $c$ and 
is minimal in the sense that there is a sequence
\begin{equation}\label{seqpost}
B\perp B_1^*\perp \cdots \perp B^*_{r(s)}\perp\bar B_s,
\quad  B\perp_c B^*_1,
\end{equation}
where either $r(s)=0$ (i.e.\,$B\perp_c\bar B_s$), 
or else all $B^*_i$, $1\le i\le r(s)$, are collinear.
\item [(ii)]
Each root $\lambda_k$ climbing over $B$ at $c$ also climbs over
 a (unique) $\bar B_s$.
\end{itemize}
\end{defn}

In Fig.1, $\{B_2,B_3\}$ is a cover of $0 \in C(B_0)$. 
In Fig.2, $\{B_2,B_3\}$ is a cover of 
$c\in C(B_0)$.
\medskip

Take a bar $\hat{B}$ of maximal height. 
Since all $\lambda_k$
are simple roots, every trunk growing on $\hat{B}$ has
bimultiplicity either $[1,0]$ or $[0,1]$.
Since $\nu_f(\hat{B})\ne 0\ne \nu_g(\hat{B})$, 
$\hat{B}$ is purely non-collinear.
It follows that every $c$ has a (unique) cover.
  
\medskip
\noindent {\bf Theorem C. }
{\em 
Let $B$ be a non-collinear bar. Take $c\in C(B)$ with
cover $\{\bar B_1,\ldots, \bar B_l\}$. Then there are exactly 
\begin{equation}\label{rootsatC}
m_B(c) + \sum_{s=1}^l [n(\bar B_s)-m(\bar B_s)] 
\end{equation}
polar roots which climb over $B$ at $c$, bounded by every $\bar
B_s$, $1\le s \le l$.
 } \medskip

Finally, let us imagine bars as tiles, 
collinear points as holes, and introduce two definitions.  
%Take $B$ with $C(B)\ne
%\emptyset$. 
By a {\it partial repair} of 
$B$,  $C(B)\ne \emptyset$, we mean a sequence,
%of postbars, 
beginning with $\tilde B_0:=B$,
\begin{equation}\label{parrep}
\tilde B_0\perp \tilde B_1 \perp \cdots \perp \tilde B_t ,
\end{equation}
where $\tilde B_t$ is purely non-collinear, $\tilde B_{i+1}$ is
supported at a collinear point of $\tilde B_i$, $0\le i\le t-1$.

The {\it repair} of $B$, denoted be $\rep (B)$, is the set of all
bars, other than $B$, which appear in any partial repair.  Let us
write
\begin{equation}\label{repair}
\rep (B) = \{\hat B_1,\ldots, \hat B_r\}.
\end{equation}

Note that if a bar appears in more than one sequence
\eqref{parrep}, it is merely taken as one element of $\rep (B)$. In
Fig.2, $\rep (B_0) = \{B_1, B_2, B_3, B_4\}$. 

\epsfxsize=9cm
\epsfysize=5cm
$$\epsfbox{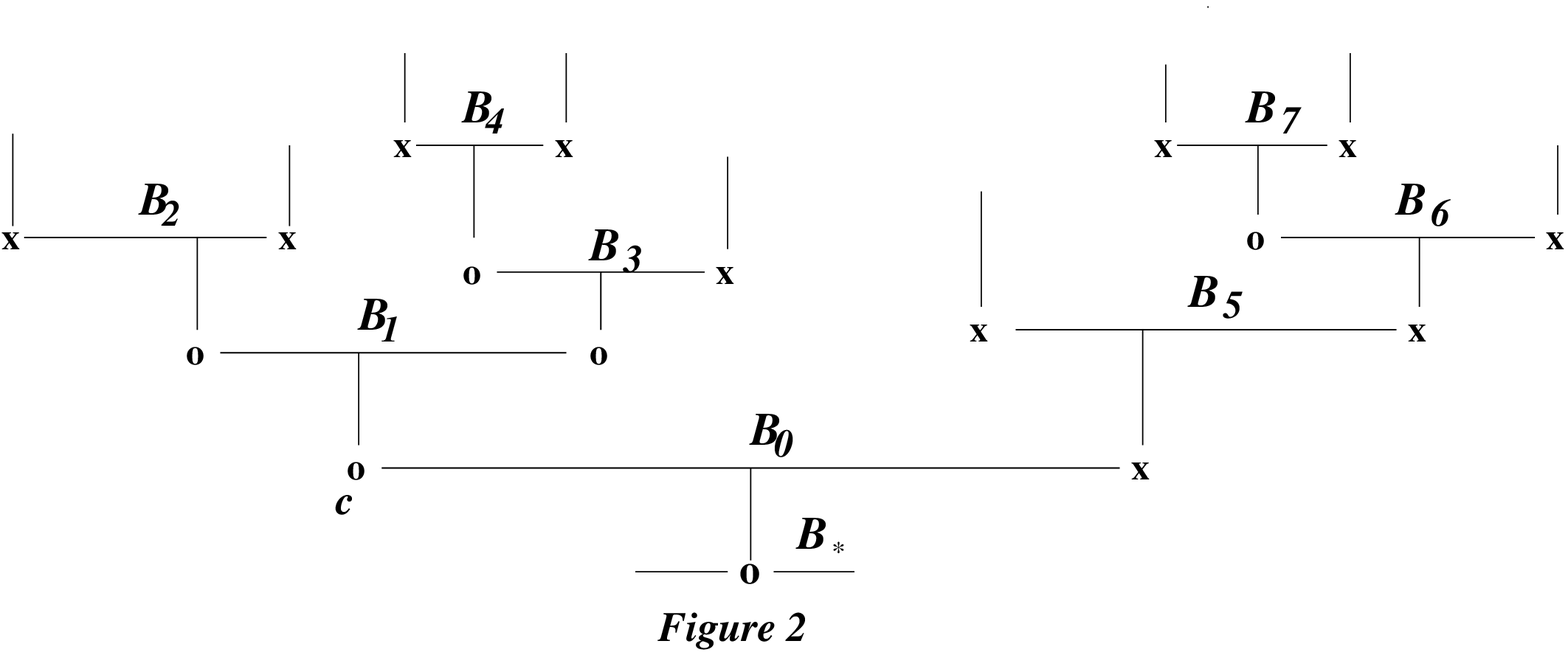}$$

A polar root is called a {\it weed} on $B$ if it climbs over $B$ at a point in 
$C(B)\cup M(B)$,
and, whenever it climbs over a bar, say $\hat B_s$, in $\rep (B)$,
it does so at a point in $C(\hat B_s) \cup M(\hat B_s)$.

Let $w(B)$ denote the number of weeds on $B$ (counting multiplicities).

\begin{cor}\label{pofC}
Let $B$ be a non-collinear bar with repair $\rep (B) = \{\hat B_1,\ldots, \hat
 B_r\}$. Then
\begin{equation}
w(B) = m(B) + \sum _{s=1}^r n(\hat B_s) .
\end{equation}
(In this formula, a collinear $\hat B_s$ yields $n(\hat B_s)=0$.) 
\end{cor}

Let $\allpost (B)$ denote the set consisting of $B$ and all 
$B'$ of finite height
lying above $B$.  We say $B'\in \allpost (B)$ is {\it basic} if
either $B'=B$ or else $B'$ is a non-collinear bar, supported at a
non-collinear point.   In Fig.2, $B_0, B_5$, and $B_6$ are the
basics in $\allpost (B_0)$.

If $B$ is purely non-collinear we define $\rep (B) := \emptyset$,
$w(B) := m(B)$.

\begin{cor}\label{total}
Let $B$ be a non-collinear bar.  Let $\{B'_1, \ldots, B'_L\}$ be the set of 
all basics in $\allpost (B)$.  Then $\T(B)=\sum_{i=1}^L w(B'_i)$.  
\end{cor}

Note that, by Theorem N, a collinear bar is never supported at a
non-collinear point, hence we have a disjoint union: 
\begin{equation}
\allpost (B) = \bigcup_{i=1}^L \,  [\{B'_i\} \cup \rep (B'_i)].
\end{equation} 

\begin{rem}
Suppose the ground bar $B_*$ is collinear. 
Take a polar root $\gamma$. Either $\gamma$ climbs over some
non-collinear bar, or else it is bounded by all non-collinear
bars of minimal height. This is because every bar of maximal height is
non-collinear.
Let $\{\bar B_1, \cdots \bar B_s\}$ be the cover of $0\in B_*$. Let us
write $J(x,y)=y^E\cdot J_*(x,y)$, $J_*$ regular in $x$, say of order $K$.
Therefore the total number of polar roots bounded by all non-collinear bars of
minimal height is $K-\sum_{i=1}^s \T(\bar B_i)$. 
\end{rem} 
\medskip

%%%%%%%%%%%%%%%%%%%%%%%%%%%%%%%%%%%%%%%%%%%%%%%%%%%%%%%%%%%%%%%%%%%

\section{Lemmas.}
\label{Lemmas}
\medskip
Take $z_k\in N(B)\cup C(B)$, with bimultiplicity $[p_k,q_k]$. 
Suppose $B\perp B^*$ at $z_k$.  

\begin{lem} \label{thesame} 
All arcs $\xi$ climbing over $B$ at $z_k$, bounded by $B^*$,
yield a constant determinant:  
\begin{equation*} 
\biggl| \begin{matrix}
\val_f(B)             & p_k  \\
\val_g(B)   & q_k  
\end{matrix} \biggr| =
 \biggl| \begin{matrix}
\val_f(\xi)             & p_k  \\
\val_g(\xi)   & q_k  
\end{matrix} \biggr| =
\biggl| \begin{matrix}
\val_f(B^*)             & p_k  \\
\val_g(B^*)   & q_k  
\end{matrix} \biggr| .
\end{equation*}
In particular, the determinants vanish if and only if $z_k\in
C(B)$; in this case, there is a common ratio:
\begin{equation*}
[\val_f(B):\val_g(B)] = [\val_f(\xi):\val_g(\xi)] =
[\val_f(B^*):\val_g(B^*)] .
\end{equation*}
\end{lem}

The proof is short.  By assumption, $$\lambda
_{B^*}(y)=\lambda_B(y)+z_ky^{h(B)}+ \cdots,\quad  \xi
(y)-\lambda_{B^*}(y)=ay^{h(B)+e}+ \cdots,$$ where $a\ne 0$,
$h(B)<h(B)+e<h(B^*)$.

Let $\zeta (y):=\lambda _B(y)+zy^{h(B)}$, $z$ a generic number.

The number of roots $\alpha_i$ with $O(\alpha_i,\xi )=h(B)+e$ is
precisely $p_k$, and $O(\alpha_i,\zeta )=h(B)$. For any other root
$\alpha_k$, $O(\alpha_k,\zeta)=O(\alpha_k,\xi).$ Hence,
%$$\xi(y) = \lambda_B(y) + z_k y^{h(B)} + \cdots;\quad h(B^*) > O
%(\xi, \lambd%a_{B^*}) %> h(B).$$
%Let $e:= O(\xi, \lambda_{B^*}) - h(B)$. Then, clearly,
\begin{equation*}
\val_f(\xi) =\sum_{j=1}^pO(\alpha_j,\zeta)+p_k\cdot e=\val_f(B) + p_k
\cdot e.
%\qquad
%\val_g(\xi) =\val_g(B) + q_k \cdot e .
\end{equation*}

Similarly, $\val_g(\xi) =\val_g(B) + q_k \cdot e$. This completes the proof.

\medskip
\begin{cor}
Take $c\in C(B)$, with cover $\{\bar B_1,\ldots,\bar B_l\}$.  
For all arcs $\xi$ which climb over $B$ at $c$, bounded by every
$\bar B_s$, $1\le s\le l$, and for all $B_j^*$ in the sequences
\eqref{seqpost} for each $\bar B_s$, there is a constant ratio: 
\begin{equation}%\label{fibr}
[\val_f(\xi):\val_g(\xi)] = [\val_f(B):\val_g(B)] = [\val_f(B_j^*):
\val_g(B_j^*)] =  [\val_f(\bar B_s):\val_g(\bar B_s)].
\end{equation} 
\end{cor}

All $B_j^*$ are collinear.  A
recursive application of Lemma \ref{thesame} completes the proof.  

\begin{lem}\label{oos}
Take any $a\notin N(B) \cup C(B)$.  For all arcs $\xi$ climbing over
$B$ at $a$, 
\begin {equation}\label{os}
\frac {y\frac d {dy} f(\xi(y),y))} {f(\xi (y),y)}  = \val_f (B) + O(y^+), 
\qquad
\frac {y\frac d {dy} g(\xi(y),y))} {g(\xi (y),y)}= \val_g (B) + O(y^+).
\end{equation} 
\end{lem}
\medskip 

Indeed, we can write
\begin{equation*}
f(\xi(y),y) = b y^e + \cdots , \quad b\ne 0,\; e=\val_f(B),
\end{equation*}
and then the lemma follows.

%%%%%%%%%%%%%%%%%%%%%%%%%%%%%%%%%%%%%%%%%%%%%%%%%%%%%%%%%%%%

\bigskip
\section{Proofs.}
\label{MandN}
\medskip

We show Theorems T, N, and C. 
The proofs depend heavily on the classical Theorem of Rouch\'e:  
\begin{equation*}
\frac 1 {2\pi i} \int _C \frac d {dz} \log \mer (z) \, dz=N-P,\quad
\text {the {\it argument index} in $C$},
\end{equation*} 
where $N, P$ denote respectively the number of
zeros and poles in a contour $C$.  

Take a non-collinear bar $B$.  Define the following meromorphic
function 
\begin{equation*}
\mer_B(z,y) := 
 \left | \begin{array}{cc}
\val_f(B)   & \sum_{i=1}^p \frac {y^{h(B)}} {x-\alpha_i(y)} \\
\val_g(B) &  \sum_{j=1}^q \frac {y^{h(B)}} {x-\beta_j(y)}  
   \end{array} \right| ,
\end{equation*}
where we have made the substitution $x=\lambda_B(y) + z
y^{h(B)}$.  It is easy to see that
\begin{equation*}
\mer_B (z) = \mer_B(z,0) \quad  (\not\equiv 0),
\end{equation*}
whence, by Rouch\'e's Theorem,
\begin{equation*}
\frac 1 {2\pi i}\int_C \frac d {dz} \log \mer_B(z,y) \, dz =
\frac 1 {2\pi i}\int_C \frac d {dz} \log \mer_B(z) \, dz, 
\quad |y|\; \text { small }.  
\end{equation*}

Using the following identities (derived from \eqref{fandg}):
\begin{equation*}
\frac {f_x} f = \sum \frac {1} {x-\alpha_i} + \frac {u_x} u ,
\qquad
\frac {g_x} g = \sum \frac {1} {x-\beta_j} + \frac {u'_x} {u'}, 
\end{equation*}
we can write $J(x,y)$ as
\begin{equation}\label{J}
J(x,y) = y\inv fg 
 \left | \begin{array}{cc}
\frac {yf_y} f   & \frac {f_x} f  \\
\frac {yg_y} g &  \frac {g_x} g  
   \end{array} \right| = y^{-h(B)-1} fg [\mer _B(z,y) + \pmer _B
(z,y)], 
\end{equation}
where
\begin{equation*}
\pmer _B(z,y):= 
 \left | \begin{array}{cc}
\frac {yf_y} f - \val_f(B)  & y^{h(B)} \frac {f_x} f  \\
\frac {yg_y} g -\val_g(B) &  y^{h(B)} \frac {g_x} g  
   \end{array} \right| + y^{h(B)}
 \left | \begin{array}{cc}
\val_f(B)   & \frac {u_x} u  \\
\val_g(B)  &  \frac {u'_x} {u'}  
   \end{array} \right| . 
\end{equation*} 

Here $\pmer _B(z,y)$ is a meromorphic function of $z,y$ which is, 
by Lemma \ref{oos}, well-defined at a generic point $(z,0)$.  

\begin{lem} \label{pmer0}
$\displaystyle \pmer _B(z,0) = 0$ for $z \notin N(B)
\cup C(B)$.  
\end{lem}

This is not obvious. Although the second summand clearly vanishes
when $y=0$, the
second column of the first determinant may not.  

Take any $\xi$ climbing over $B$ at a point $a\notin N(B)
\cup C(B)$. Let us evaluate $\pmer_B$ along $\xi$.  The
first determinant vanishes by Lemma \ref{oos}.  Since $\xi$ is
arbitrary, Lemma \ref{pmer0} follows.

\begin{cor}
Take $a\in \C$, and a small $\varepsilon >0$. Take $y \in \C$, $|y| \ll 
\varepsilon$. Then 
\begin{eqnarray}\label{ints} 
& \frac 1 {2\pi i}\int_{|z-a|=\varepsilon}
\frac d {dz} \log [\mer_B(z,y) + \pmer_B(z,y)] \, dz 
\\
& = \frac 1 {2\pi i}\int_{|z-a|=\varepsilon}
\frac d {dz} \log \mer_B(z) \, dz = \mu_B(a)\notag.
\end{eqnarray}
\end{cor}
\medskip

Now let us take $a=z_k$ on $B$. There are $\tau_B(z_k)$ roots of 
$f(x,y)g(x,y)$ in the contour $|z-z_k|=\varepsilon$. It follows
from the above corollary that $J(x,y)$ has
$\tau_B(z_k)+\mu_B(z_k)$ roots in the contour. Since
$\varepsilon$ is arbitrarily small, these roots must all climb over
$B$ at $z_k$. This completes the proof of Theorem T.

\medskip
To prove Theorem N, we can permute the indices in \eqref{allzk}, if necessary,
so that $z=z_1$ with bimultiplicity $[p_1,q_1]$. Let $B^*$ be
the postbar of $B$ supported at $z_1$, and let
\begin{equation*}
N(B^*)= \{z_1^*, \ldots, z_s^*\}, \quad
C(B^*)= \{z_{s+1}^*, \ldots, z_{s+t}^*\}, 
\end{equation*}
where $s=n(B^*)$, $t=c(B^*)$.  Then, clearly, 
\begin{equation*}
\sum _{k=1}^{s+t} p^*_k = p_1, \qquad
\sum _{k=1}^{s+t} q^*_k = q_1, 
\end{equation*}
and, since $z_1$ is non-collinear, 
\begin{equation} \label{D*}
D^*:= \sum_{k=1}^{s+t} \,  \biggl| \begin{matrix}
\val_f(B) & p_k^*  \\
\val_g(B) & q_k^* 
\end{matrix} \biggr| =
\biggl| \begin{matrix}
\val_f(B) & p_1  \\
\val_g(B) & q_1
\end{matrix} \biggr|
\ne 0 . 
\end{equation}
Moreover, 
\begin{equation*}
D^*=   \biggl| \begin{matrix}
\val_f(B^*) & p_1  \\
\val_g(B^*) & q_1 
\end{matrix} \biggr| = \sum_{k=1}^s \, 
\biggl| \begin{matrix}
\val_f(B^*) & p_k^*  \\
\val_g(B^*) & q_k^* 
\end{matrix} \biggr|; 
\end{equation*}
the first equality follows from Lemma \ref{thesame},
the second holds because $z_{s+j}^*$ are collinear.
\medskip

Now consider $\mer_{B^*} (z)$, the rational function associated
to $B^*$.  Its numerator is a polynomial of degree $s-1$ with
leading coefficient $D^*\ne 0$. We have proved \eqref{m+1=n}.

The second part of Theorem N follows from \eqref{m+1=n} and
Theorem T.  Indeed, let us calculate the number of polar roots
climbing over $B$ at $z$ and the same number over $B^*$.  
By Theorem T
the former equals $p_1+q_1-1$
and the latter equals $\sum  p^*_k + \sum q^*_k + m(B^*) -n(B^*)$.

For Theorem C the ideas of the proof are the same as before.  Calculating the
argument index of  $\mer_B = y^{h(B)+1} f\inv g\inv J(x,y) $ on $B$ 
at $c\in C(B)$ yields
$m_B(c)$.  Calculating the index on each $\bar B_s$, $1\le s\le
l$, in a large contour yields $m(\bar B_s) - n(\bar B_s)$.  The
total deficit is the number in \eqref{rootsatC}, proving Theorem
C.

For Corollary \ref {pofC}, we take a cover of each $c\in C(B)$.  
If $\bar B_s$ in \eqref{seqpost} is not purely non-collinear we take
the covers of its collinear points.  This process is repeated
recursively.  Now, adding all the collinear bars $B_j^*$
appearing in \eqref{seqpost} to the above covers yields the
repair $\rep (B)$.  Hence Corollary  \ref {pofC} follows from
Theorems C and T.

Take a polar root $\gamma$ which climbs over $B$. By Theorem N, 
there is a {\it unique} basic $B'_s$ in
$\allpost (B)$ for which $\gamma$ is accounted for in $w(B'_s)$,
$1\le s\le L$.  This proves Corollary \ref {total}.

%%%%%%%%%%%%%%%%%%%%%%%%%%%%%%%%%%%%%%%%%%%%%%%%%%%%%%%%%%%%%%
\bigskip
\section{What Theorem C Does Not Say}
\label{C}
\medskip 

Theorem C does not say precisely where the polar roots leave the
tree. The number of polar roots given in Theorem C is
determined by the tree, but their orders of contact with the tree need
{\it not} be.  It follows that the contact structure of the two curve
germs, {\it i.e.} $T(f,g)$, does not give full information on how to factor 
the Jacobian into irreducible factors in $\C\{x,y\}$. (See Theorem F below.) 
We shall use two examples to show that in our case, contrary to the one 
function case see Section \ref{Onefunction}, the way these polar
 roots leave the tree need not be an invariant of the tree; 
the coefficients of the $\lambda_i$'s may also play a r\^ole. 

First, consider Example \ref{ex1.1}. The tree model is shown in Fig.1,
with $\val_f (B_2) = \val_g(B_2) = \val_f (B_3) = \val_g(B_3) = E+e+3$, $B_1$ 
being collinear.  By Theorem T, there are four polar roots climbing
over $B_0$, all at $0$,

In the following, we take $e<E<2e$.

To decide how many polar roots climb over $B_1$, we put $x=zy^{e+1}$. 
An easy calculation yields 
$$
\frac{yg_y} {g}-\frac{yf_y} {f}
=2y^{E-e} \Bigl [ \frac {ezy^{2e-E}}
{z^2y^{2e}-1}-\frac{(E-e)A(z-1)} {(z-1)^2-A^2y^{2(E-e)}}-\frac
{(E-e)B(z+1)} {(z+1)^2-B^2y^{2(E-e)}}\Bigr ];
$$ 
and 
$$\frac {g_z}
{g}-\frac {f_z} {f}=2y^{E-e}\Bigr [ \frac {y^{2e-E}}
{z^2y^{2e}-1}+\frac {A} {(z-1)^2-A^2y^{2(E-e)}}+\frac {B}
{(z+1)^2-B^2y^{2(E-e)}}\Bigl ] . 
$$ 
Hence 
$$J(x,y)=y^{-e-2}\cdot f\cdot 
g\cdot 
 \left | \begin{array}{cc}
 \frac {yf_y} {f}   & \frac {f_z} {f} \\
\frac {yg_y} {g }&   \frac {g_z} {g}
   \end{array} \right|
=2y^{E-2e-2}\cdot f\cdot g\cdot \Delta(z,y),
$$ 
where 
$$\Delta (z,0)=
 \left | \begin{array}{cc}
{2e+3}   & \frac 1 {z-1}+ \frac 1 {z+1} \\ 
\frac {(e-E)A} {z-1}+\frac {(e-E)B} {z+1} 
&  \frac {A} {(z-1)^2}+\frac {B} {(z+1)^2}
   \end{array} \right|$$
$$=(z^2-1)^{-2}[(A+B)(2E+3)z^2+2(A-B)(E+e+3)z+(2e+3)(A+B)].$$

Observe that if $A+B\ne 0$, there are two zeros.  
This means that
two polar roots climb over $B_1$, the remaining two are bounded
by $B_1$.
If, however, $A+B=0$, then there is only one zero. This means
that one polar root climbs over $B_1$, 
three are bounded by $B_1$.

Thus, in general, 
{\it one cannot tell the positions of polar roots 
relative to collinear bars. }

\begin{example}\label{ex6.1}
Take integers $e$, $E$, $N$, $2e>E>e>0$, $N\ge 0$. Let $$
f(x,y):=[x^2-y^{2(e+1)}][(x-y)^2-y^{2(e+1+N)}],\qquad
g(x,y):=[x+y^{E+1}][x+y].$$ The tree has four bars, $B_*$, $B_1$,
$B_2$, $B_3$; $h(B_1)=1, h(B_2)=e+1, h(B_3)=e+1+N$; and
$$
\mer_{B_1}(z)=\frac8{z^2-1},\qquad
\mer_{B_2}(z)=\frac{-2(e+2)}{z(z^2-1)}, $$ 
where $B_2$ is supported at a
collinear point $0\in B_1$.

By Theorem C, three polar roots, say
$\gamma_i$, $1\le i\le 3$, climb over $B_1$ at $0$, bounded by $B_2$.

Let us write $x=Xy$.  The arcs 
$\eta_i(y):=y^{-1}\gamma_i(y)$, $i=1, 2, 3$, are Newton-Puiseux
roots of the equation 
$$
 X^3 (8 + \cdots) - Xy^E (2(E+2) + \cdots)- y^{2e}(2(e+2) + \cdots)=0 ,
$$ 
where the doted terms are in the maximal ideal.  If $3E<4e$ then 
the Newton Polygon of this equation has vertices $(3,0)$, $(1,E)$, and
$(0,2e)$.  Then  two $\eta_i$'s
have order $\frac E2$, one has order $2e-E$.
Thus, two polar roots have order $\frac
E2+1$, one has order $2e-E+1$.

We can take $e=7$. Then $E_1:=8$, $E_2:=9$ both satisfy the
above inequality. Let $$ g_k(x,y)=(x+y^{E_k+1})(x+y), \qquad k=1,2.$$
Then $T(f,g_1)=T(f,g_2)$, but, as $E_1\ne E_2$, the polar roots
split away from the trees at different heights between $B_1$ and
$B_2$.
\end{example}

The above examples show that the factorization can depend on  the
coefficients of the Newton-Puiseux roots, not merely on the contact structure.

%%%%%%%%%%%%%%%%%%%%%%%%%%%%%%%%%%%%%%%%%%%%%%%%%%%%%%%%%%%%%%%%%%
\bigskip
\section{Factors of $J(x,y)$ in $\C\{x,y\}.$}
\label{F}
\medskip
We can introduce an additional structure on $T(f,g)$ and use it to define a 
factorisation of $J(x,y)$ in $\C\{x,y\}$. The
factors are not irreducible, in general.

\begin{defn}%\label
Take $B$, $\bar B$. We say $B$ is {\it conjugate} to $\bar B$, writing as $B
\thicksim \bar B$, if, and only if $h(B)=h(\bar B)$ and there exists an
irreducible $p(x,y)\in \C\{x,y\}$, of which one (Newton-Puiseux) root
climbs over $B$ and one climbs over $\bar B$.
\end{defn}

\begin{lem} \label{conjugate} 
Suppose $B \thicksim \bar B.$ Take any irreducible $q(x,y)\in
\C\{x,y\}.$ If $q(x,y)$ has a root climbing over $B$ then it also has
a root climbing over $\bar B$.
\end{lem}
Proof. Take an integer $D$ such that the roots of $p(x,y)$ and $q(x,y)$ can
all be written in the following form:
\begin{equation}\label{lambda}
\lambda (y):=c_1y^{\frac {n_1}{D}} +c_2y^{\frac{n_2}{D}} + \cdots,\qquad
0<\frac {n_1}{D}<\frac{n_2}{D}<\cdots.
\end{equation}

Of course, here we allow $D$, $n_1$, $n_2$, \dots, to have common
factors.

Let $\theta$ be any $Dth$ root of unity: $\theta^D=1$. Each $\theta$
yields a transformation (conjugation) on arcs of form
$\eqref{lambda}$:
$$\theta(\lambda)(y):=c_1\theta^{n_1}y^{\frac{n_1}{D}}
+c_2\theta^{n_2}y^{\frac{n_2}{D}} +\cdots.$$

As in \cite {walker} (p.107), it is easy to see that the $\theta$'s permute
transitively the
roots of $p(x,y)$, and also that of $q(x,y)$. The contact order is preserved; 
in particular,
$$O(\alpha,\beta)=O(\theta(\alpha),\theta(\beta)),\quad
\text{if} \quad p(\alpha(y),y)=q(\beta(y),y)=0.$$

Lemma \ref{conjugate} follows immediately. That $\thicksim$ is an
equivalence relation also follows.
\medskip

Thus, one can simply use any irreducible component of $f(x,y)g(x,y)$
as $q(x,y)$ to identify an equivalence class of bars at any given height.
\medskip

Let $\mathfrak{B}:=\mathfrak{B}(f,g)$ denote the set of all
equivalence classes of bars.
\medskip

Take $\mathbb B\in \mathfrak B$. If some $B\in \mathbb B$ is collinear
(resp. non-collinear) then every $\bar B \in \mathbb B$ is collinear
(resp. non-collinear); in this case we say $\mathbb B$ is collinear
(resp. non-collinear).

Let $C(\mathfrak B)$, $N(\mathfrak B)$ denote respectively the
collinear and non-collinear classes of bars.

Let us take an integer $D$ so that $\lambda_k$, $1 \le k \le N$, can
all be written in the form \eqref{lambda}.
\medskip

Take $\theta$, $\theta^D=1$. Take $z\in B$ (Convention \ref{copy}).
 If $h(B)=\frac{n}{D}$, then $\theta(z)=\theta^nz$, {\it i.e.},
$$\theta(\lambda_B(y)+zy^{h(B)})=\lambda_{\bar
  B}(y)+\theta^nzy^{h(\bar B)}; \quad \bar B\thicksim B.$$

If $z_k\in C(B)$ (resp.$N(B)$), having bimultiplicity $[p_k,q_k]$,
then $\bar z_k:=\theta(z_k)\in C(\bar B)$ (resp.$N(\bar B)$), having
bimultiplicity $[\bar p_k,\bar q_k]=[p_k,q_k]$. Hence
$$ \nu_f(B)=\nu_f(\bar B), \quad \nu_g(B)=\nu_g(\bar B).$$

Observe also that 
$$\mer_{\bar B}(\theta^nz)=\theta^{-n}\mer_B(z),$$
whence $\theta$ induces a bijection between the pure mero-zeros of $B$
and $\bar B$, preserving the mero-multiplicity.

Now take $\mathbb B$, non-collinear, and consider the product
$$
P_{\mathbb B}(x,y):= \prod_j[x-\gamma_j(y)],
$$
taking over all $j$ such that $\gamma_j$ leaves the tree on some 
$B\in \mathbb B$.
\medskip

\begin{lem} \label{lfactor}  
Take any non-collinear $\mathbb B\in \mathfrak B$. Then 
$P_{\mathbb B}(x,y)\in \C\{x,y\}$, and hence
$$P(x,y):=\prod_{\mathbb B \in N(\mathfrak B)}P_{\mathbb B}(x,y)$$
is a factor of $J_{(f,g)}(x,y)$ in $\C\{x,y\}$.
\end{lem}

Take a polar root $\gamma$, leaving the tree on $B\in\mathbb B$. Take
an irreducible $p(x,y)\in\C\{x,y\}$ having $\gamma$ as a root.  
Every root of $p(x,y)$ leaves the tree on some $\bar B\in \mathbb B$. Hence 
$P_{\mathbb B}(x,y)$ is a product of factors like $p(x,y)$. This completes 
the proof.
\medskip

Take $\mathbb B$, non-collinear, $h(\mathbb B):=h(B)>0$, $B \in
\mathbb B.$  Consider the product
$$ Q_{\mathbb B}(x,y):=\prod [x-\gamma_j(y)],$$
taking over all $j$ such that $\gamma_j$ climbs over some $B\in \mathbb
B$ at some $c\in C(B)$, and is bounded by every bar of the cover of
$c$. If $\mathbb B$ is purely non-collinear, we define $Q_{\mathbb B}(x,y):=1.$

As for the ground bar $B_*$, it may be collinear or not. Let
$$ Q_{B_*}(x,y):=\prod [x-\gamma_j(y)],$$
taking over all $j$ such that $\gamma_j$ is bounded by every
non-collinear bar of minimal height.
\medskip

\noindent{\bf Theorem F.}
{\em
The above defined $Q_{\mathbb B}(x,y)$ and $Q_{B_*}(x,y)$ are in $\C
\{x,y \}$, and hence
\begin{equation}\label{factor} 
J_{(f,g)}(x,y)=unit \cdot y^E \cdot Q_{B_*}(x,y) \cdot \prod_{\mathbb B \in
  N(\mathfrak B)}P_{\mathbb
  B}(x,y) \cdot Q_{\mathbb B}(x,y),\quad E\ge 0,
\end{equation}
$h(\mathbb B)>0$, is a decomposition in $\C\{x,y\}$. }
\bigskip

%%%%%%%%%%%%%%%%%%%%%%%%%%%%%%%%%%%%%%%%%%%%%%%%%%%%%%%%%%%%%%%%%%%%%%%%%%%%%%

\section{Tree-invariants and Mero-invariants}\label{Inv}
\medskip

In this section we analyse how the factorisation \eqref{factor} is determined 
by the tree.   For this we first define what is meant by a tree-invariant of 
the pair $(f,g)$. 
We say $(f,g)$ and $(f',g')$ are \emph{equivalent} if the following two
 conditions are satisfied.
\medskip

\noindent{\bf Condition 1.}
{\em
(Compare \cite {Zariski}.) The roots of $f$ and $f'$, and of $g$ and $g'$, are
in one-one correspondence: $\alpha_i \post \alpha'_i$, $1\le i \le p$;
$\beta_j \to \beta'_j$, $1\le j \le
q$. The contact order is preserved:
$$O(\alpha_i, \alpha_k)=O(\alpha'_i,\alpha'_k),\quad 
O(\alpha_i,\beta_j)=O(\alpha'_i,\beta'_j),\quad 
O(\beta_j,\beta_s,)\,=O(\beta'_j,\beta'_s).$$}

This is equivalent to saying that there is a bijection between the bars
and trunks of $T(f,g)$ and that of $T(f',g')$ which preserves the
heights and bimultiplicities.  If a bar $B$ corresponds to a bar $B'$, then 
$\val_f (B) = \val _{f'} (B')$ and $\val_g (B) = \val _{g'} (B')$ and 
hence $B$ is collinear iff so is $B'$.  

In particular, if $B$ corresponds to $B'$, both non-collinear, then $C(B)\cup 
N(B)$ and
$C(B')\cup N(B')$ are in one-one correspondence, say $z_k\post z'_k$,
as in (2.1), such that 
$$ \Delta_B(z_k)=\Delta_{B'}(z'_k),\qquad 1\le k\le l.$$ 
\noindent{\bf Condition 2.}
{\em
If $B$ corresponds to $B'$, non-collinear, then
$m(B)=m(B')$; and if $c\in C(B)$ corresponds to $c'\in C(B')$ then $
m_B(c)=m_{B'}(c')$. }
\medskip

We say that an object associated to $(f,g)$ is a \emph{tree-invariant of 
$(f,g)$} if it
depends only on the equivalence class of $(f,g)$. 

\begin{rem}\label{pure} 
If $T$ is a trunk of bimultiplicity $[s,0]$ or $[0,t]$. Let $B$ be the bar on 
top of $T$.  
Then Condition 2 for $B$ is satisfied automatically by Corollary \ref{ST}.  
\end{rem}

Take $B\in\mathbb B$, non-collinear. Let $m^*(B)$ denote the number of
pure mero-zeros on $B$ (counting multiplicities). Let 
$$m^*(\mathbb B):=\sum_{B\in \mathbb B}m^*(B);\quad \nu_f(\mathbb
B):=\nu_f(B);\quad \nu_g(\mathbb B):=\nu_g(B).$$

Let $P_{\mathbb B}$ denote the curve germ defined by $P_{\mathbb
B}(x,y)=0$.

As a consequence of Condition 1, an equivalence between $T(f,g)$ and 
$T(f',g')$ induces a bijection
between $\mathfrak B(f,g)$ and $\mathfrak B(f',g')$. Suppose $\mathbb B$
corresponds to $\mathbb B'$. Then $\mathbb B$ is non-collinear if, and only
if $\mathbb B'$ is; in this case $m^*(\mathbb B)=m^*(\mathbb B')$, by 
Condition 2. 
\medskip

\noindent{\bf Theorem I.}
{\em  
If $\mathbb B$ is non-collinear, then 
$$ I(C_f,P_{\mathbb B})=\nu_f(\mathbb B)m^*(\mathbb B);\qquad  
I(C_g,P_{\mathbb B})=\nu_g(\mathbb B)m^*(\mathbb B);$$
and 
$$
I(C,P_{\mathbb B})=[\nu_f(\mathbb B)+\nu_g(\mathbb B)]m^*(\mathbb B), 
$$
where $C_f = f\inv (0)$, $C_g = g\inv (0)$, and $C=C_f\cup C_g$.  
These intersection multiplicities are tree-invariants of $(f,g)$. 
Similarly the orders (in $x$) of the factors in \eqref{factor} are 
tree-invariants. } 
\smallskip 

If a polar root leaves the tree on a non-collinear bar, it must do so at a pure
 mero-zero (Corollary \ref{P}). Hence the intersection multiplicities 
$I(C_f,P_{\mathbb B}), I(C_g,P_{\mathbb B})$ are tree-invariants. The above 
formulae follow from Corollary \ref{M}.

By Theorem C, the number of the $\gamma_j$'s in the definition of
$Q_{\mathbb B}(x,y)$ is a tree-invariant. This number is the order. 
The last statement of Theorem I is proved. 

\medskip
\noindent{\bf Attention.}
{\em
 Example \ref{ex1.1}, as analysed in \S \ref{C}, shows that similar
 intersection multiplicities defined by $Q_{\mathbb B}(x,y)$ need
  not be tree-invariants. We have no formula for them either.  }
\medskip

\noindent{\bf Addendum.}
{\em  Take a non-collinear $\mathbb B$. Let 
$$P_{\mathbb B}(x,y)=p_1(x,y)^{e_1} \cdots p_s(x,y)^{e_s}, \quad e_i\ge 1,$$
be the irreducible decomposition of $P_{\mathbb B}(x,y)$ in $\C\{x,y\}.$
Then, clearly,
$$s \le e_1+ \cdots +e_s \le m^*(B), \quad  B\in \mathbb B.$$   }

Note that $m^*(B)$ is a tree-invariant, but $s$, $e_i$ are
not.

The following example shows that the Zariski equisingularity types (\cite {Zariski}) of 
$P_{\mathbb B}(x,y)$ need not be 
tree-invariants.

\begin{example}\label{ex8.2}(Compare Example(4.4) in \cite{kuo-lu}.)
Take $g=g'=g''=y$, and
$$f=x^3-y^4, \quad f'=x^3-y^4-3xy^5, \quad f''=x^3-y^4-3xy^6.$$
Note that $(f,g)$, $(f',g')$, and $(f'',g'')$ are equivalent. There is
only one bar of height $\frac {4}{3}$ in each tree, for which
$$ P_{\mathbb B}(x,y)=x^2;\quad P_{{\mathbb B}'}(x,y)=x^2-y^5; 
\quad P_{{\mathbb B}''}(x,y)=x^2-y^6.$$
\end{example}
\medskip

We now define the {\it truncation} relative to $T(f,g)$, 
$P_{\mathbb B}^\top$, as follows.
Take any arc $\xi$, not a root of $f\cdot g$, leaving the
tree on $B$ at $a$. We call
$$ \xi^{\top}(y):=\lambda_B(y)+ay^{h(B)}$$
the {\it truncation} of $\xi$ {\it relative} to $T(f,g)$. 
If $\xi$ is one of the roots $\lambda_k$, we define $\xi^{\top}:=\xi$.

Take $h(x,y)\in\C\{x,y\}$, regular in $x$, say of order $r$, with factorisation
$$h(x,y)=unit\cdot \prod_{i=1}^r[x-\xi_i],\quad O(\xi_i)>0.$$

We define the truncation of $h(x,y)$ {\it relative} to $T(f,g)$ by
$$h^{\top}(x,y):=\prod_{i=1}^r[x-\xi_i^{\top}(y)].$$

A simple conjugation argument shows that $h^{\top}(x,y)\in
\C\{x,y\}$.
Note that in Example \ref{ex8.2}, $P_{\mathbb B}^{\top}=P_{{\mathbb
    B}'}^{\top}=P_{{\mathbb B}''}^{\top}=x^2.$
\medskip
 
We say $(f,g)$ and $(f',g')$ are {\it mero-equivalent} if the
following condition is also satisfied.

\noindent{\bf Condition 3.}
{\em
Suppose $B$ corresponds to $B'$, non-collinear. Then there
exists a bijection between the pure mero-zeros of $B$ and
$B'$, preserving the mero-multiplicity.}
\medskip

It is easy to see that in this case, $P_{\mathbb B}^{\top}$ and
$P_{{\mathbb B}'}^{\top}$ are  Zariski equisingular.

\medskip
We say a function germ associated to $T(f,g)$ is a {\it mero-invariant} if
its Zariski equisingularity type depends only on the mero-equivalence
class of $(f,g)$.
\medskip

\noindent{\bf Addendum.}
{\em
In Theorem I the intersection multiplicities do not change when 
$P_{\mathbb B}$ is replaced by the curve germ defined by $P_{\mathbb B}^{\top}
(x,y)=0$. The truncations $P_{\mathbb B}^{\top}(x,y)$ are mero-invariants.  }

%%%%%%%%%%%%%%%%%%%%%%%%%%%%%%%%%%%%%%%%%%%%%%%%%%%%%%%%%%%%%%%%%%%%%%%%
\bigskip
\section{Discussions}\label{Disc}
\medskip

\subsection{} 
Let $f,g:(\C^2,0)\post (\C,0)$ be given. We can take a generic
constant $c$ and substitute $y$ by $y+cx$. Then $f, g$ are
mini-regular in $x$ (\cite {kuo-lu}), that is, $E_1=E_2 =0$ and 
$O(f)=p$, $O(g)=q$.  
In the factorisations
(\ref{fandg}),
$O(\alpha_i)\ge 1$ and $O(\beta_j)\ge 1$.
Since $c$ is generic, $J(x,y)$ is also mini-regular in $x$, say of
order $m$. Hence there are exactly $m$ polar roots, all of order $\ge
1$. These are called the ``generic'' polar roots.

The following example shows in particular that $m$ need not be a
tree-invariant.

\begin{example}\label{ex9.1}
Take $f(x,y)=x^2-G(x,y)^2$, $g(x,y)=x-2G(x,y)$, where
$$G(x,y):=\int_0^y(x-t^2)^2dt, \qquad G(x,0)=0, \qquad G(0,y)=\frac15
y^5.$$ Then $J(f,g)=-2(2x-G)G_y=-2(2x-G)(x-y^2)^2,$ having $m=3.$ 

The ground bar is collinear. The other bar of $T(f,g)$ is non-collinear, having
height 5. Two of the three polar roots are bounded by this bar; they
leave the tree at height 2.

Next, take $f'=x^2-G(x^2,y)^2$, $g'=x- 2G(x^2,y)$. Then $T(f',g')=
T(f,g)$, and $$ J(f',g')=-2(2x-G(x^2,y))(x^2-y^2)^2, \quad m=5.$$
Four polar roots are bounded by the bar of height 5; they
leave the tree at height 1.
\end{example}

\smallskip
\subsection{One function case} \label{Onefunction}

 Our results generalise that of the one function case. Let $f(x,y)$
be regular in $x$, $g(x,y)=y$. Then
$J(x,y)=f_x$, $T(f,g)=T(f)$, the tree-model defined in
(\cite {kuo-lu}). All bars are purely non-collinear (Corollary
\ref{ST}). Theorem T reduces to Lemma(3.3) of \cite {kuo-lu}: if
there are $l$ trunks growing on $B$, then there are $l-1$ polar roots
leaving the tree on $B$.

Let us first suppose $f(x,y)$ is irreducible (Merle's case). Let $\{\frac
{q_1}{p_1},\frac{q_2}{p_1p_2},\cdots,\frac{q_g}{p_1\cdots
  p_g}\}$ denote the Puiseux characteristic sequence of $f$, as in
\cite {kuo-lu} (p.308). The tree $T(f)$ has a simple form (\cite {kuo-lu}). 
Since $f$ is irreducible, all bars of the same height are obviously
conjugate. Let ${\mathbb B}_s$ denote the class of
bars with height $\frac{q_s}{p_1\cdots p_s}$, $1\le s\le g$. 

Take ${\mathbb B}_s$. Take any polar root, $\gamma$, which leaves the
tree on some $B\in {\mathbb B}_s$. Note that $ \nu_f(B)=O(f(\gamma
(y),y)$.
There is a sequence of postbars, $B_1\perp B_2\perp \cdots \perp B_s$,
where $B_s:=B$, with heights $h(B_k)=\frac {q_k}{p_1\cdots p_s}$, $1\le
k\le s$.

Given $k$, there are $p_k$ trunks growing on $B_k$, and $p_k-1$
polar roots leaving the tree on $B_k$.

Hence the total number of such $\gamma$ is $(p_s-1)p_{s-1}\cdots
p_1$. It follows that
$$f_x=unit\cdot \prod_{i=1}^gP_{{\mathbb B}_i}(x,y);\quad I(C_f,P_{{\mathbb
    B}_s})=\nu_f({\mathbb B}_s)[p_s-1]p_{s-1}\cdots p_1.$$
This is Merle's theorem (\cite {merle}) (true also for non-generic polars).

Next consider the general case where $f$ may be reducible
(Garcia-Barroso's case). Let $H$ denote the set of contact orders,
$O(\alpha_i,\alpha_j)$, between the Newton-Puiseux roots $\alpha_i$ of $f$. 
Note 
that $H$ is just the set of heights $h(B)$, for all $B$.

Take $h\in H$. Unlike Merle's case, the bars with the same height $h$ may
have more than one conjugate classes. Let them be denoted by
${{\mathbb B}_1^{(h)}},\cdots, {{\mathbb B}_{r(h)}^{(h)}}.$ Theorem F implies:
$$f_x=unit\cdot \prod_{h\in H}\prod_{i=1}^{r(h)}P_{{\mathbb
    B}_i^{(h)}}(x,y).$$
This is Garcia-Barroso's theorem (\cite {Barro}) (true also for
non-generic polars). An equivalence class of bars is a "black point"
of Eggers (\cite {egg}), the height is the "valuation".

Example \ref{ex8.2} also exposes the Pham phenomenon in the one
function case. Note that $f$, $f'$, and $f''$ all have Puiseux
exponent $\frac43$, hence $T(f)=T(f')=T(f'')$. There are two polar
roots in each tree. For $f$, they generate the same factor
$x$: $f_x=3x^2$; for $f'$, an irreducible factor: $f'_x=3(x^2-y^5)$;
and for $f''$, two distinct factors: $f''_x=3(x-y^3)(x+y^3)$.

This explains why in the factorisations of Merle, Garcia-Barroso and
our Theorem F, the factors need not be invariants of the tree.

A brief summary:\,In the one function case, the factors of $f_x$ in
$\C \{x,y\}$ are not invariants of $T(f)$ (\cite{Pham}), but the set of
contact orders, $C(f,f_x)$, {\it is} (\cite{kuo-lu}). For a general
pair $(f,g)$, however, if there are collinear points or bars, then the
contact order set need not be an invariant (\S \ref {C}).  

\medskip
\subsection{Meromorphic Case} 
%%%%%%%%%%%%%%%%%%%%%%%%%%%%%%%%%%%%%%%%%%%%%%%%%%%%%%%%%%%%

One can easily extend our results to the meromorphic case to generalise the 
results of Assi, see \cite{assi} (and also \cite{abhyankar-assi} or 
\cite{abhyankar-assi2}).

Let $F(X,Y)$, $G(X,Y)$ be a given pair of monic polynomials 
in $X$: 
$$F(X,Y)=X^p+a_1(Y)X^{p-1}+ \cdots +a_p(Y),\; G(X,Y)=X^q+ \cdots +b_q(Y),$$
where $a_i(Y), b_j(Y)\in \C(Y)^*$, the field of fractional 
power series of $Y$.  The Newton-Puiseux Theorem asserts that $\C(Y)^*$ is 
algebraically closed (\cite{walker}, p.98).

For simplicity, let us assume $a_i$, $b_j$ are Laurent series.

Take a large integer $s$. We then use the 
substitution $X=xy^{-s}$, $Y=y$, and define 
$$
f(x,y):=F(xy^{-s},y),\quad g(x,y):=G(xy^{-s},y). 
$$
Then, clearly, $y^{ps}f(x,y)$ and $y^{qs}g(x,y)$ are holomorphic and
hence $f$ and $g$ can be factored in the form (0.1) with $u = y^{-ps}$, 
$u' = y^{-qs}$, (not units). Moreover,
$$
XYJ_{(F,G)}(X,Y)=xyJ_{(f,g)}(x,y);
$$ 
and $X=\xi (Y)$ is a root of $J(X,Y)$ if and only if $x=y^s\xi (y)$
is one of $J(x,y)$.

We have thus reduced the meromorphic case to the holomorphic case. As before, 
we can define the tree-model, $T(f,g)$, the associated rational functions, 
etc..  An important observation is that in this case, $\nu_f(B)$, $\nu_g(B)$ 
can be negative, or zero. As an example, let us take 
$$
F(X,Y)=X^4 -Y^{-2}X^2 +1,\quad G(X,Y)=X^2 - 
Y^{-1}X,$$ 
with $X=xy^{-2}$, $Y=y$. The tree $T(f,g)$ has a bar, $B$, with 
$h(B)=2$, $\nu_f(B)=\nu_g(B)=0$. A bar $B$ with $\nu_f(B)=\nu_g(B)=0$ is 
obviously collinear; thus, if $B$ is also of maximal height, then any 
collinear point which lies below $B$ has no cover.

Theorems T, N, and C remain true.  Corollary(2.7), Corollary(2.8) and 
formula (2.14) remain true if we assume that for every bar, $B'$, lying above 
$B$, $\nu_f(B')\ne 0\ne \nu_g (B')$.

Theorem(JF1), in \cite {abhyankar-assi}, 
(and the Theorem in [2], Section 9, p9,) yields factors under the hypothesis 
$\Omega_B(G)=1$  (the notation of \cite 
{abhyankar-assi}).  The contact set (defined in \cite {abhyankar-assi}, 
p.129, lines 5 and 21) is
essentially the tree-model defined in \cite {kuo-lu}, a bud (\cite 
{abhyankar-assi}) is essentially a bar in \cite {kuo-lu};
Theorem(DF1) of \cite{abhyankar-assi} is Lemma(3.3)  of \cite{kuo-lu}.  
In our terminology, the condition $\Omega_B(G)=1$ means 
that there is no Newton-Puiseux root of $G$ climbing over $B$, i.e.\;$B$ has
bimultiplicity $[s,0]$, so this case is not different to the one function 
case.   In particular, (see Corollary \ref{ST},) 
$B$ is purely non-collinear and $\tau(B) >0$, except (in the meromorphic 
case) when $\nu_B(G)=0$, that is $S(G,B)=0$ in the terminology of
\cite{abhyankar-assi}.    In general, however, the tree $T(F\cdot G)$ of the 
product function can obviously have many ``buds'' $B$  
which do not have this property.

%%%%%%%%%%%%%%%%%%%%%%%%%%%%%%%%%%%%%%%%%%%%%%%%%%%%%%%%%%%%%%%%%%%%%%%%%%

\end{document}